\documentclass[12pt]{amsart}
\pdfoutput=1
\usepackage{amssymb}
\usepackage{amsmath}
\usepackage{amsfonts}
\usepackage[usenames]{color}
\usepackage{graphicx}
\usepackage{array}
\usepackage{psfrag}
\usepackage{color}
\usepackage{ulem}

\makeatletter
 
 \@addtoreset{equation}{section}
\makeatother

\newtheorem{theorem}{Theorem}

\theoremstyle{definition}

\theoremstyle{remark}

\numberwithin{equation}{section}

\newcommand{\bibtitle}[1]{\textit{#1}}

\begin{document}

\title{Hadwiger numbers of self-complementary graphs}

\date{\today}
\author{Andrei Pavelescu}
\address{University of South Alabama, Mobile, AL 36688}
\email{andreipavelescu@southalabama.edu}

\author{Elena Pavelescu}
\address{University of South Alabama, Mobile, AL 36688}
\email{elenapavelescu@southalabama.edu}

\begin{abstract}
The Hadwiger number of a graph $G$, denoted by $h(G)$, is the order of the largest complete minor of $G$. 
 A graph is said to be self-complementary if it is isomorphic to its complement.
We prove that for all $n\equiv 0,1 (\text{mod 4})$ and for all $ \lfloor \dfrac{n+1}{2} \rfloor \le h \le \lfloor \dfrac{3n}{5}\rfloor $, there exists a self-complementary graph $G$ with $n$ vertices whose Hadwiger number is $h$.\end{abstract}
\vspace{0.1in}
\maketitle

\section{Introduction}

A minor of a simple undirected graph $G$ is a graph $H$ that can be obtained from $G$ through a series of vertex deletions, edge deletions, and edge contractions. 
The Hadwiger number of a simple undirected graph $G$, denoted by $h(G)$, is the order of the largest complete minor of $G$. 
The Hadwiger conjecture, one of the famous problems in graph theory, states that for any graph $G$,  $\chi(G) \le h(G)$, where $\chi(G)$ is the chromatic number of $G$. While the validity of the conjecture is still unknown in general, Girse and Gillman \cite{GG} and Rao and Sahoo \cite{RS} proved the conjecture to be true for self-complementary graphs. 
A graph is said to be self-complementary if it is isomorphic to its complement.
 Nordhaus and Gaddum \cite{NG} proved that for a self-complementary graph $G$ on $n$ vertices, $\chi(G)\le \lfloor \dfrac{n+1}{2} \rfloor$.
 Rao and Sahoo \cite{RS} showed that $ \lfloor \frac{n+1}{2} \rfloor \le h(G)$.
The same lower bound for $h(G)$ was independently found by the authors in \cite{PP}.
In \cite{RS} and \cite{PP},  different classes of examples show that this lower bound is attained.

Motivated by a conjecture of Kostochka \cite{K}, Stiebitz \cite{St} proved that if $G$ is a graph on $n$ vertices, then $h(G)+h(cG)\le \lfloor \frac{6n}{5} \rfloor$, where $cG$ is the complement of $G$.
This implies that for a self-complementary graph with $n$ vertices $G$,  $h(G) \le \lfloor \frac{3n}{5}\rfloor $. Given these upper and lower bounds for the Hadwiger numbers of self-complementary graphs,  Rao and Sahoo \cite{RS} asked whether each integer within this allowable range is realized as a Hadwiger number of a self-complementary graph.
In what follows, we answer this question in the positive.
We first show that the $\lfloor \frac{3n}{5}\rfloor$ upper bound is realized for all $n\equiv 0,1 (\text{mod 4})$, thus proving the upper bound for the Hadwiger number is sharp.
Then we use induction to prove the following:

\begin{theorem} For all $n\equiv 0,1 (\text{mod 4})$ and for all $ \lfloor \dfrac{n+1}{2} \rfloor \le h \le \lfloor \dfrac{3n}{5}\rfloor $, there exists a self-complementary graph $G$ with $n$ vertices whose Hadwiger number is $h$.
\label{theorem1}
\end{theorem}

We prove Theorem \ref{theorem1} in the next two sections.

\section{Constructions for the upper bound}

Self-complementary graphs with $n$ vertices exist only for $n\equiv 0,1 (\text{mod 4}) $. One is prompted by the upper bound of $\lfloor \frac{3n}{5}\rfloor $ to consider the remainders of $n$ modulo 5. The Chinese Remainder Theorem yields ten cases that need to be considered. The following construction, introduced in \cite{RS}, provides examples of self-complementary graphs with maximum Hadwiger number in eight of the cases.
 
\begin{figure}[htpb!]
\begin{center}
\begin{picture}(250, 100)
\put(0,0){\includegraphics[width=3.5in]{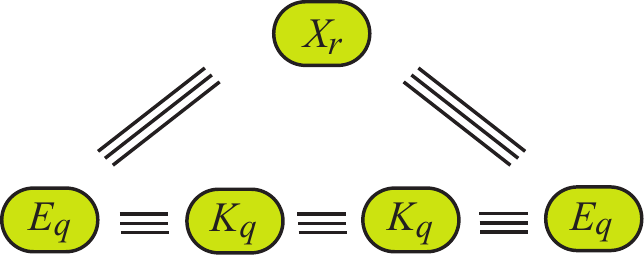}}
\end{picture}
\caption{\small A self-complementary graph with $n = 4q+r$ vertices.}
\label{Fig_5cycle}
\end{center}
\end{figure}

Consider $X_r$ a self-complementary graph on $r$ vertices, $K_q$ the complete graph on $q$ vertices, and $E_q$ the graph with $q$ vertices and no edges.
The graph $G$ with $n= 4q+r$ vertices described in Figure \ref{Fig_5cycle} is self-complementary.
In this figure, the triple lines mark that all edges between the respective subgraphs are present (i.e. a complete bipartite graph).

By assigning different values to $r$ and $q$ in the graph described in Figure \ref{Fig_5cycle}, we obtain self-complementary graphs of maximal Hadwiger number  for all $n\equiv 0,1 (\text{mod 4})$, except $n\equiv 12,17 (\text{mod 20} ).$
These graphs are presented in Table \ref{Table_SCregular}.
We explain why each of these eight graphs attain the maximum Hadwiger number of $\lfloor\frac{3n}{5}\rfloor$.

\begin{table}[htpb!]
\centering
\caption{}
\label{Table_SCregular}
\begin{tabular}{|l|l|l|l|}
\hline
 $n=4q+r$ & $r$ & $q$ & $h=\lfloor\frac{3n}{5}\rfloor$   \\ \hline
 $20s$& $4s$ & $4s$ & $12s$  \\ \hline
 $20s+1$& $4s+1$  & $4s$ &$12s$    \\ \hline
$20s+4$ &$4s$  &$4s+1$  & $12s+2$  \\ \hline
$20s+5$ & $4s+1$  & $4s+1$ & $12s+3$   \\ \hline
$20s+8$ &  $4s$& $4s+2$ & $12s+4$   \\ \hline
$20s+9$ & $4s+1$  & $4s+2$ & $12s+5$  \\ \hline
$20s+13$ &  $4s+1$& $4s+3$ & $12s+7$  \\ \hline
$20s+16$ &  $4s+4$& $4s+3$ & $12s+9$  \\ \hline
\end{tabular}
\end{table}

Let $p$ denote the minimum between $q$ and $r$. For each copy of $E_q$, contract $p$ disjoint edges between $p$ of its vertices and the same $p$ vertices of $X_r$. This way, we obtain a $K_{2q+p}$ minor of $G$. Notice that for all eight graphs in Table \ref{Table_SCregular}, $2q+p$ equals the upper bound $\lfloor\frac{3n}{5}\rfloor$, thus $h(G)=\lfloor\frac{3n}{5}\rfloor$.

For $n= 20s +12$ and $n=20s+17$, we found that no values of $r$ and $q$  yield graphs whose Hadwiger number is $\lfloor\frac{3n}{5}\rfloor $ ($12s+7$ and $12s+10$, respectively).
For these two cases, we provide different classes of examples.\\

\begin{figure}[htpb!]
\begin{center}
\begin{picture}(200,180)
\put(0,0){\includegraphics[width=2.5in]{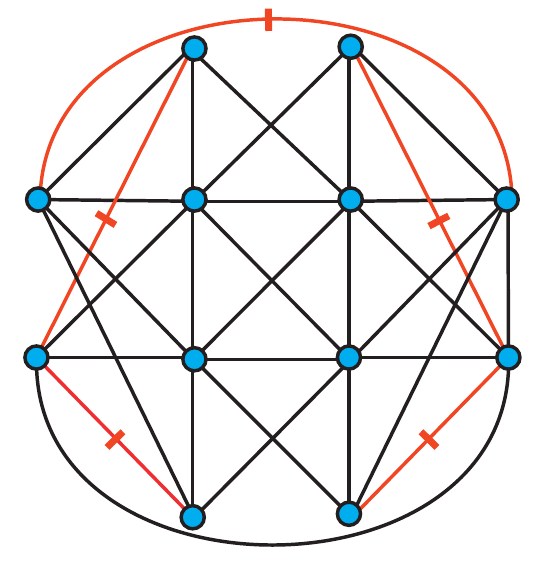}}
\end{picture}
\caption{\small A self-complementary graph with $n = 12$ vertices containing a $K_7$ minor. This minor can be obtained by contracting the five marked edges.}
\label{Fig_SC12K7}
\end{center}
\end{figure}

For $n=12$, $\lfloor \frac{3n}{5} \rfloor=7$. The self-complementary graph with 12 vertices presented in Figure \ref{Fig_SC12K7} contains a $K_7$ minor obtained by contracting the five marked edges.


\begin{figure}[htpb!]
\begin{center}
\begin{picture}(400,185)
\put(0,0){\includegraphics[width=6in]{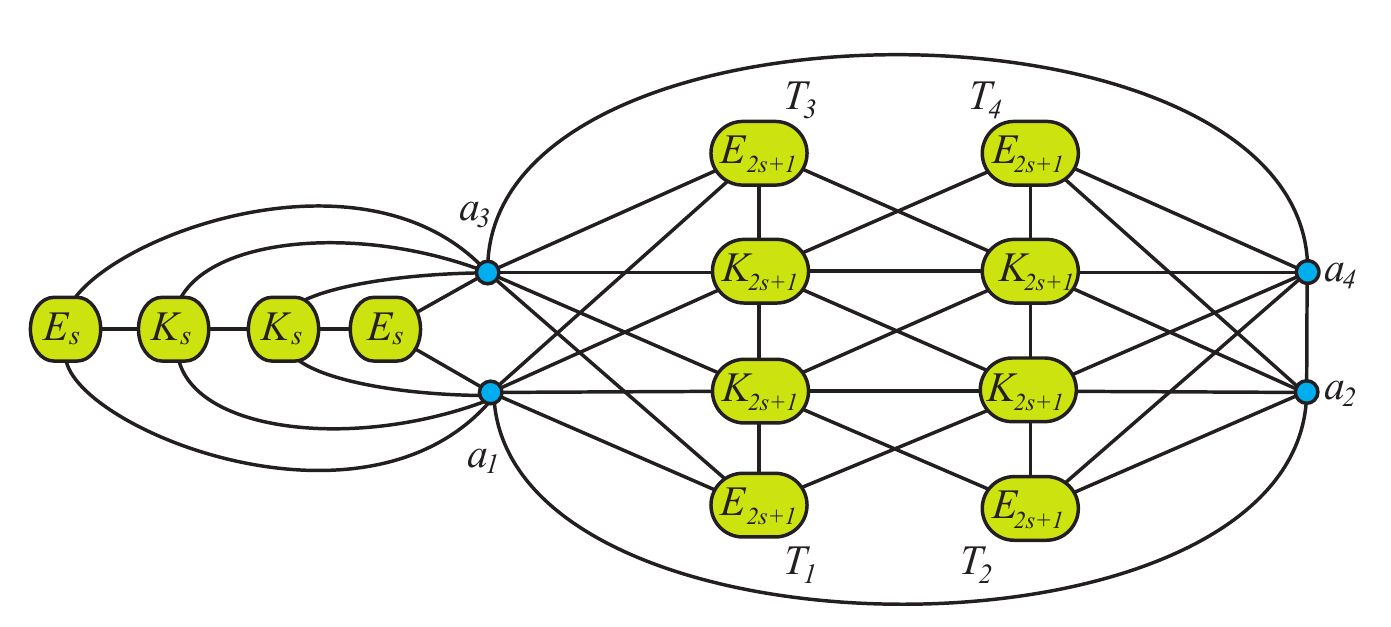}}
\end{picture}
\caption{\small A self-complementary graph on $n = 20s+12$ vertices containing a $K_{12s+7}$ minor.  
All the drawn edges represent complete bipartite graphs, and $a_1, a_2, a_3$, and $a_4$ are single vertices.
All edges between $T_1, T_2, T_3, T_4$ and the two copies of $E_s$ and the two copies of $K_s$ are also included in the graph. These are not drawn to preserve readability. }
\label{Fig_SCbig}
\end{center}
\end{figure}

For $n=20s+12$ with $s\ge 1$, let $G$ be the self-complementary graph obtained from the graph in Figure \ref{Fig_SCbig}, together with all edges between the vertices of $T_1, T_2, T_3, T_4$ and the vertices of the two copies of $E_s$, and all edges between the vertices of $T_1, T_2, T_3, T_4$ and the vertices of the two copies of $K_s$.
For a natural number $m$, $K_m$ denotes the complete graph with $m$ vertices and $E_m$ the graph on $m$ vertices with an empty edge set.  
This graphs admits a $K_{12s+7}$ minor obtained by performing the following sequence of edge contractions:
\begin{itemize}
\item For $i\in \{1,3\}$, choose one vertex $t_i$ of $T_i$ and contract $2s$ disjoint edges between the remaining vertices of $T_i$ and the $2s$ vertices of the two copies of $K_s$;
\item For $i\in \{2,4\}$, choose one vertex $t_i$ of $T_i$  and contract $2s$ disjoint edges between the remaining vertices of $T_i$ and the $2s$ vertices of the two copies of $E_s$;
\item the edge $a_3a_4$ contracts to form a single vertex $a_3=a_4$;
\item the edges $a_1t_1$ and $a_1t_3$ contract to form the vertex $a_1=t_1=t_3$.
\item the edges $a_2t_2$ and $a_2t_4$ contract to form the vertex $a_2=t_2=t_4$;
\end{itemize}

Note that the four copies of $K_{2s+1}$ induce a $K_{8s+4}$ subgraph of $G$. 
The contractions of edges in between the vertices of the $T_i$'s and those of the $K_s$' and those of $E_s$', respectively, take place within bipartite graphs.
 As such, there are many possible choices of the set of contracted edges. 
 This sequence of edge contractions produces a $(8s+4)+4s+3 =12s+7$ complete minor.

\begin{figure}[htpb!]
\begin{center}
\begin{picture}(250, 95)
\put(0,0){\includegraphics[width=3.5in]{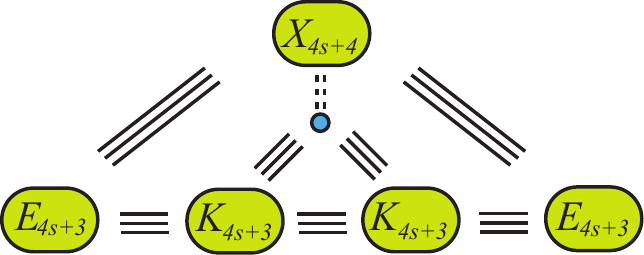}}
\put(123, 40){$a$}
\end{picture}
\caption{\small A self-complementary graph on $n = 20s +17$ vertices. The double dashed line marks that the highlighted vertex is adjacent to exactly half the vertices in $X_{4s+4}$.}
\label{Fig_SC20s17K12s10}
\end{center}
\end{figure}

For $n= 20s +17$, $s\ge 0$, we build a self-complementary graph $G$ on $n$ vertices by adding a vertex $a$ to the graph in Figure \ref{Fig_5cycle} with $r=4s+4$ and $q=4s+3$. See Figure \ref{Fig_SC20s17K12s10}.


 Consider the subgraph of $G$ induced by the vertices of $X_{4s+4}$.
Since the average degree of the vertices of $X_{4s+4}$ is $\frac{4s+3}{2}$, exactly half of them have degree less than $\frac{4s+3}{2}$ (see \cite{PP}). Let $b$ be one of these vertices. The vertex $a$ neighbors all the vertices of the two copies of $K_{4s+3}$, and half of the vertices of $X_{4s+4}$, namely those of small degree, $b$ included.
To obtain the $K_{12s+10}$ complete minor, contract $4s+3$ disjoint edges between one copy of $E_{4s+3}$ and the $4s+3$ vertices of $X_{4s+4}$ except $b$ , another  $4s+3$ disjoint edges between the other copy of $E_{4s+3}$ and the same $4s+3$ vertices of $X_{4s+4}$, and the edge $ab$.


\section {Inductive Step} Let $n\ge 12$, $n\equiv 0, 1$(mod 4). Assume that for each integer $k$ between $\lfloor \frac{n+1}{2}\big \rfloor$ and $\lfloor \frac{3n}{5} \rfloor$, there exists a self-complementary graph $G$ on $n$ vertices whose Hadwiger number is $k$. Using the construction in Figure \ref{Fig_5cycle} with $X_r=G$ and $q=1$, we obtain a self-complementary graph on $n+4$ vertices and Hadwiger number $k+2$. As 
\[\bigg \lfloor \frac{n+4+1}{2} \bigg \rfloor-\bigg \lfloor \frac{n+1}{2}\bigg \rfloor=2\,\,\, 
 \text{and} \,\,\, 
 \bigg \lfloor \frac{3(n+4)}{5} \bigg \rfloor-\bigg \lfloor \frac{3n}{5}\bigg \rfloor\le 3,\]  it follows that, for every $k$ in between $ \lfloor \frac{n+4+1}{2} \rfloor$ and  $\lfloor \frac{3(n+4)}{5} \rfloor$, one inductively builds an example of a self-complementary graph on $n+4$ vertices and Hadwiger number $k$, with the exception of $k=\lfloor \frac{3(n+4)}{5} \rfloor$ when $\lfloor \frac{3(n+4)}{5} \rfloor-\lfloor \frac{3n}{5}\rfloor=3$. 
 Since the case of the upper bound $k=\lfloor \frac{3(n+4)}{5} \rfloor$ was already covered in the previous section, the proof is complete by induction on $n$. 
\bibliographystyle{amsplain}

\end{document}